\documentclass[ba]{imsart}

\RequirePackage[authoryear]{natbib}%% uncomment this for author-year citations
\RequirePackage[colorlinks,citecolor=blue,urlcolor=blue]{hyperref}
\RequirePackage{graphicx}

\usepackage{amsthm}
\usepackage{amsmath}
\usepackage[utf8]{inputenc} % allow utf-8 input
\usepackage[T1]{fontenc}    % use 8-bit T1 fonts
\usepackage{hyperref}       % hyperlinks
\usepackage{url}            % simple URL typesetting
\usepackage{booktabs}       % professional-quality tables
\usepackage{amsfonts}       % blackboard math symbols
\usepackage{nicefrac}       % compact symbols for 1/2, etc.
\usepackage{microtype}      % microtypography
\usepackage{lipsum}		% Can be removed after putting your text content
\usepackage{graphicx}
\usepackage{natbib}
\usepackage{soul}
\usepackage{color}

\newcommand\medidaAleatoria{\mathcal{P}}

\newcommand\Amostra{\mathbf{X}_n}
\newcommand\posterior{\Pi_{\theta|\Amostra}}

\newcommand\prior{\Pi_0}

\newcommand\particao{\mathcal{B}}

\newcommand\ArvoreDePolya{PT(\particao, \mathcal{A})}

\newcommand\eje{E_j^*}

\newcommand\e{\mathbb{E}}

\newcommand\len{l(\epsilon)}
\newcommand\nne{N_\epsilon}
\newcommand\nze{N_{\epsilon0}}
\newcommand\nnu{N_{\epsilon1}}
\newcommand\nnel{N_{\epsilon'}}
\newcommand\nzel{N_{\epsilon'0}}
\newcommand\nnul{N_{\epsilon'1}}

\newcommand\al{a_{\len}}

\newcommand\BO{B_{\epsilon0}}
\newcommand\BU{B_{\epsilon1}}
\newcommand\BE{B_\epsilon}
\newcommand\eps{\epsilon}
\newcommand\epsl{\epsilon'}
\newcommand\epss{\epsilon^{*}}

\newcommand\estimador{\hat{H}(\Amostra)}

\newcommand\ex{\mathbb{E}_0}

\newcommand\parcelaR{\mu(\nze, \nnu)}
\newcommand\parcelaRl{\mu(\nzel, \nnul)}

\newcommand\fbo{P_X(\BO)}
\newcommand\fbu{P_X(\BU)}
\newcommand\fbe{P_X(\BE)}

\newtheorem{assumption}{Assumption}

\newtheorem{prop}{Proposition}
\newtheorem{theorem}{Theorem}
\newtheorem{corollary}{Corollary}
\newtheorem{lemma}{Lemma}
\newtheorem{example}{Example}
\newtheorem*{theorem*}{Theorem}

\begin{document}
\begin{frontmatter}
\title{Kullback-Leibler Consistency of $p$-dimensional Pólya Tree Posteriors and differential entropy estimation}
%\title{A sample article title with some additional note\thanksref{t1}}
\runtitle{Pólya Tree Posteriors: KL Consistency differential entropy estimation}
%\thankstext{T1}{A sample additional note to the title.}

\begin{aug}
%%%%%%%%%%%%%%%%%%%%%%%%%%%%%%%%%%%%%%%%%%%%%%%
%% Only one address is permitted per author. %%
%% Only division, organization and e-mail is %%
%% included in the address.                  %%
%% Additional information can be included in %%
\author[A]{\fnms{Fernando}~\snm{Corrêa}\ead[label=e1]{fptcorrea@gmail.com}\orcid{0000-0002-8050-8554}},
\author[B]{\fnms{Rafael}~\snm{Bassi Stern}\ead[label=e2]{rbstern@gmail.com}\orcid{0000-0002-4323-4515}}
\and
\author[C]{\fnms{Julio}~\snm{Michael Stern}\ead[label=e3]{jmstern@gmail.com}\orcid{0000-0003-2720-3871}}
%%%%%%%%%%%%%%%%%%%%%%%%%%%%%%%%%%%%%%%%%%%%%%
%% Addresses                                %%
%%%%%%%%%%%%%%%%%%%%%%%%%%%%%%%%%%%%%%%%%%%%%%
\address[A]{Institute of Mathematics and Statistics, University of São Paulo\printead[presep={,\ }]{e1}}
\address[B]{Institute of Mathematics and Statistics, University of São Paulo\printead[presep={,\ }]{e2}}
\address[C]{Institute of Mathematics and Statistics, University of São Paulo\printead[presep={,\ }]{e3}}
\runauthor{Corrêa, Stern and Stern}
\end{aug}

\begin{abstract}
We exploit the multiplicative structure of Pólya Tree priors to establish novel consistency results on $p$-dimensional trees, conditions to obtain Kullback-Leibler minimax contraction rates for univariate density estimation and a representation theorem of entropy functionals of Pólya Tree posteriors. These results motivate a novel differential entropy estimator that is consistent under mild conditions on large dimensions.
\end{abstract}



\end{frontmatter}
%%%%%%%%%%%%%%%%%%%%%%%%%%%%%%%%%%%%%%%%%%%%%%
%% Please use \tableofcontents for articles %%
%% with 50 pages and more                   %%
%%%%%%%%%%%%%%%%%%%%%%%%%%%%%%%%%%%%%%%%%%%%%%
%\tableofcontents

\section{Introduction}

Pólya Trees \citep{mauldin_polya_1992} are widely studied stochastic processes that draw random probability measures. They are convenient prior distributions in nonparametric Bayesian inference since they are conjugate, mathematically tractable, and allow the modeling of both absolutely continuous and non-continuous distribution functions. Pólya Tree mixtures and generalizations were applied in many different scenarios such as hypothesis testing, survival analysis, and directional statistics \citep{lavine_aspects_1992,kraft_class_1964,ferguson_prior_1974, lima1, castillo-spike-an-slab}.

%We present a novel posterior consistency result for Pólya Trees and a novel differential entropy estimator motivated by. Both results are connected by an 

One reason for theoretical interest in Pólya Trees is posterior consistency, 
a major goal in Bayesian nonparametric inference. There are many definitions of consistency, and they materialize the concentration of the posterior around the true data generating process. Usual theorems on the consistency of Pólya Trees explore particular cases of Schwarz's theorem, or the application of the more recent multi-scale approach \citep{barron_consistency_1999,castillo_bernsteinvon_2015}.

The usual Bayesian nonparametric setting for density estimation considers $\Amostra = (X_1, ..., X_n)$ an i.i.d. sample distributed according to a random density $\theta$ over a sample space $\mathcal{X}$.  $\theta$ is sampled by a measure $\prior$ with support on the set of densities supported on $\mathcal{X}$. For a measurable set $S$ the posterior measure $\posterior$ is then given by

\begin{equation}\label{eq:definition-of-posterior}
\posterior(S) = \frac{\int_S \prod_{i=1}^n \theta(x_i) d\prior(\theta)}{\int\prod_{i=1}^n \theta(x_i) d\prior(\theta)}\tag{posterior measure}
\end{equation}

In this setting, let $f_0$ be the true value of $\theta$. A posterior $\posterior$ is \textbf{consistent} with respect to some semimetric $d$ if

\begin{equation}\label{eq:cons}
\posterior\left(\lbrace\theta: d(f_0, \theta) > \epsilon\rbrace\right) \rightarrow 0\tag{posterior consistency condition}
\end{equation}

in $f_0$ probability or almost surely with respect to an infinite sample of $X_i$. The posterior $\posterior$ is \textbf{strongly consistent} if the convergence is almost sure and \textbf{weakly consistent}, otherwise. Furthermore, a sequence $\epsilon_n$ is a \textbf{posterior contraction rate} if 

\begin{equation}\label{eq:cons2}
\posterior\left(\lbrace\theta: d(f_0, \theta) > \epsilon_nM_n\rbrace\right) \rightarrow 0\tag{contraction rate}
\end{equation}

for all sequences $M_n \rightarrow \infty$. Contraction rates are further classified as weak or strong whether this convergence occurs in probability or almost surely.

The theoretical machinery mobilized to obtain consistency results depends on $d$. Consequently, this choice defines the conditions that $\prior$ must satisfy. When $d$ metricizes weak convergence, Schwarz's Theorem requires $f_0$ to be on the Kullback-Leibler support of $\prior$ in order to achieve consistency. For other choices of $d$, much stronger requirements are needed. For the Hellinger norm $d_H$ or the total variation norm $d_{TV}$, existing results require $\prior$ to assign small probability to rough densities. Applying the recent multi-scale approach ensures that similar conditions also imply convergence with respect to the supremum norm $d_{\infty}$ \citep{walker_new_bayesian, Ghosal_vanderVaart_2017, barron_consistency_1999}.

In this paper, we present a novel technique for obtaining consistency results on Pólya Tree priors. It allows us to expand existent consistency results and to determine sufficient conditions under which the posteriors contracts at minimax rate for the $L_1$ and Kullback-Leibler semimetric. Let $K(f, \theta)$ denote the Kullback-Leibler divergence between a density $f$ and $\theta$ sampled from a Pólya Tree. Our main result characterizes the asymptotic behavior of the posterior expectation of $K(f,\theta)$. For clarity, we will postpone the precise statement to the following sections.

\textbf{Theorem 1.}\label{teo:prob-cons}
    \textit{If $f_0$ is a density on $\mathcal{X}$ that belongs to the Kullback-Leibler support of a Pólya Tree $\prior$ and $\theta$ is almost surely bounded away from $0$ and $\infty$, then the posterior expectation $\e_{\theta \sim \posterior}[K(f_0, \theta)] \rightarrow 0$ both almost surely and in $L_1$. Furthermore if $\mathcal{X} = [0,1]$ and $f_0$ is $\alpha$-Holder and bounded away from $0$ and infinity there are Pólya tree priors such that}
    $$\e_{\Amostra \sim f_0}\left[\e_{\theta \sim \posterior}[K(f_0, \theta)]\right] = O\left(n^{-\frac{2\alpha}{1+2\alpha}}\right).$$

Theorem 1 yields posterior consistency, since $\eps^{-1}\e_{\theta \sim \posterior}[K(f_0, \theta)]$ is an upper bound for $\posterior(K(f_0, \theta) > \eps)$. Strong consistency is achieved as this quantity vanishes almost surely and in $L_1$. Moreover, Theorem 1 establishes the minimax rate of convergence for $\e_{\theta \sim\posterior}[K(f_0, \theta)]$. Furthermore, we demonstrate throughout the paper that our condition on $\prior$ is less restrictive than those previously considered. The conditions on $f_0$ are standard in this context.

This paper also establishes a new consistent estimator for
the differential entropy, $H(f_0) = -\int f_0(t) \log f_0(t) dt$. This is a challenging problem since, for example, Schwartz's Theorem does not ensure the posterior convergence of $H(\theta)$ to $H(f_0)$. Indeed, consistent Bayesian differential entropy estimation has only been proposed recently \citep{entropy_bjps,castillo_bernsteinvon_2015}. We propose the following estimator:
\begin{align*}
    \estimador =- \mathbb{E}_{\theta\sim\posterior}\left[\sum_{i=1}^n \frac{\log \theta(X_i)}{n}\right].
\end{align*}

\textbf{Theorem 2.}
    \textit{Let $d\geq1$, $f_0$  be a density on $[0,1]^d$ with finite differential entropy and $\prior$ be a Pólya Tree that samples densities on $[0,1]^d$. Under regularity conditions on $\prior$ and $f_0$,
    $\estimador$ converges to $H(f_0)$ in probability.} %%Furthermore, the posterior distribution of $\sum_{i=1}^n \frac{\log \theta(X_i)}{n} | \Amostra$ is approximately Gaussian.}

Most differential entropy estimators fall in large classes of estimators: nearest neighborhood  \citep{KozachenkoLeonenko1987}, plug-in \citep{SricharanWeiHero2013}, sample spacings \citep{Vasicek1976}, and histograms \citep{HallMorton1993}. Much of the literature discusses traditional estimators, particularly the Kozachenko-Leonenko statistic. This approach achieves consistency in probability provided that $f_0$ has finite differential entropy and additional tail conditions on $f_0$. Also, entropy estimation in general dimension $d$ is challenging. Our estimator provides a novel approach for differential entropy estimation that attains consistency in probability across any dimension with mild conditions on $f_0$.

Theorems 1 and 2 are derived from a novel representation of $K(f,\theta)$ as a random series.

\textbf{Theorem 3.}\textit{ If $f$ belongs to the weak Kullback-Leibler support of a Pólya Tree $\prior$ that sample densities almost surely bounded away from $0$ and $\infty$ , then $K(f, \theta)$ satisfies}

    $$K(f, \theta) = \sum_{i=1}^\infty P_i$$

\textit{$\prior$-almost surely where $P_i$ are mutually independent random variables. Also the sequence $\e_{\theta \sim \posterior}[K(f,\theta)]$ is a supermartingale with respect to the natural filtration.}

This representation allows new conclusions about Polya Trees.
%It establishes that $\theta$ being bounded away from $0$ and $\infty$ is a property observed in a large class of Pólya Tree priors.
Its analytical tractability enables the computation of posterior credible regions based on the quantity $K(f,\theta)$.
The representation also establishes that $E_{\theta \sim \posterior}[K(f,\theta)]$ is a supermartingale, which is a desirable property for nonparametric priors.

%vOur main idea is noting that the posterior of $\frac{\sum \log \theta(x_i)}{n}$ might concentrate very slowly around $\int \freal(t) \log \freal(t)dt$. We shall show that is quantity is approximately Gaussian under the posterior. Therefore its distribution might be characterized by its mean denoted $\frac{\eps}{n}$ and its variance $\frac{\vps}{n^2}$. The speed of concentration is related to the smoothness of $\theta$. For a certain class of rough priors on $\theta$:

% $$\frac{\left(\frac{\sum \log \theta(x_i)}{n} - \int \freal(t) \log \freal(t)dt\right)}{\frac{\sqrt{\vps}}{n}}\rightarrow -\infty$$ 

% and this leads to inconsistency. This follow from the fact that  $ \int \freal(t) \log \freal(t)dt$ might be estimated applying the Law of Large Numbers to $\sum \log \freal(X_i)$, which decrease rate is $\sqrt{n}$. 

% Through our arguments we also note that smooth $\medidaAleatoria$ estimate $\int \freal(t) \log \freal(t)dt$ from below, which is convenient for our hypothesis test.

The paper is organized as follows. In Section \ref{sec:priors-notation} we present notation and preliminary results that will be referenced throughout the paper. In Section \ref{sec:main-results} we establish our main results. We begin by presenting Theorem \ref{teo:prob-cons} that powers up all of our arguments. Then we proceed to the discussions regarding the proof of Theorem \ref{teo:prob-cons} in Section \ref{sec:KL-inc}. Finally, we analyze the properties of the differential entropy estimator in Section \ref{sec:properties-estimator}. We present summarized proofs throughout the paper and full proofs are available in the appendix.

\section{Basic definitions}
\label{sec:priors-notation}

\subsection{Pólya Trees}

In this section, we establish the notation and review classic results on Pólya Trees that are used throughout the paper. We follow the notation of \cite{ghosh_bayesian_2003} and restate results from \cite{Ghosal_vanderVaart_2017}. Throughout the paper, we adopt \( \theta \) as the symbol for a random density sampled by a Pólya Tree. This choice emphasizes that the parameter of interest is the density, rather than the random measure sampled from the Polya Tree. This focus is feasible only when the latter is almost surely absolutely continuous with respect to the Lebesgue measure $\lambda$. In this section we establish conditions that ensure that $\theta$ is well-defined.

The following notation is used:
\begin{itemize}
    \item $E_j^*$: set of all binary sequences of length $j$. We refer to elements $\epsilon = (\epsilon_1 \epsilon_2 \ldots \epsilon_j) \in E_j^*$ as sequence of binary digits $\epsilon_i \in \{0,1\}$, $1 \leq i \leq j$. For example, $(\epsilon_1 \epsilon_2 \epsilon_3) = (101) \in E^*_3$ is a binary sequence of length 3.
    \item $E_j = \bigcup_{i=1}^j E_j^*$: set of binary sequences of length less or equal than $j$.
    \item $E = \{\emptyset\}\bigcup\left(\bigcup_{i=1}^\infty E_j^*\right)$: set of all finite-length binary sequences including the empty sequence.
    \item $l(\epsilon)$: length of a binary sequence $\epsilon \in E$.
    \item $\mathcal{X}$: a compact subset of $\mathbb{R}^d$, $d\geq1$.
    \item $B(\mathcal{X})$: the Borel $\sigma$-algebra of $\mathcal{X}$.
    \item $\lambda$: the Lebesgue measure on $(\mathcal{X}, B(\mathcal{X}))$.
    \item $\mathcal{B} = \{B_{\emptyset}\}\cup\{P_1, P_2, \ldots\}$: a collection of partitions $P_j = \{B_\epsilon: \epsilon \in E^*_j\}, j \in \mathbb{N}$, such that
    \begin{enumerate}
    \item $B_{\emptyset} = \mathcal{X}$;
    \item $\{B_{\epsilon0}, B_{\epsilon1}\}$ is a partition of $B_\epsilon$ for all $\epsilon \in E$ and
    \item $\mathcal{B}$ generates $B(\mathcal{X})$.
    \end{enumerate}
    \item $\mathcal{A}$: a set of positive real numbers indexed by $E$
    $$\mathcal{A} = \{\alpha_\epsilon, \epsilon \in E\}, \ \alpha_{\epsilon} \in \mathbb{R}_+.$$
\end{itemize}  

A \textbf{Pólya Tree} prior is a probability measure $\Pi$ parametrized by $\mathcal{B}$ and $\mathcal{A}$ that generates a random probability measure $\medidaAleatoria$ over $\mathcal{X}$ which satisfies:

\begin{enumerate}
    \item All random variables in the set $\{\medidaAleatoria(B_{\epsilon0}|B_\epsilon): \epsilon \in E\}$ are independent, and
    \item For all $\epsilon \in E$, $\medidaAleatoria(B_{\epsilon0}|B_\epsilon) \sim \text{Beta}(\alpha_{\epsilon0}, \alpha_{\epsilon1}).$
\end{enumerate}

A random measure $\medidaAleatoria$ drawn according to the measure $\Pi$ is denoted as $\medidaAleatoria \sim \Pi$, or equivalently, $\medidaAleatoria \sim \ArvoreDePolya$. Furthermore, we represent the beta-distributed random variables that constitute $\medidaAleatoria$ by
$$Y_\epsilon = Y_{\epsilon_1 \ldots \epsilon_k} := \medidaAleatoria(B_{\epsilon_1 \ldots \epsilon_k}|B_{\epsilon_1 \ldots \epsilon_{k-1}}).$$

%\end{definition}

Pólya Trees are conjugate priors. Let $\Amostra|\medidaAleatoria$ be an i.i.d. sample of observations distributed accordingly to $\medidaAleatoria$. %Let $\amostra$ be an observation from $\Amostra$.
If $\medidaAleatoria \sim \ArvoreDePolya$ then the posterior $\medidaAleatoria|\Amostra$ is also sampled by a Pólya tree $PT(\mathcal{B}, \mathcal{A}_{\Amostra})$ where
$$\mathcal{A}_{\Amostra} = \{\alpha_{\eps} + N_\epsilon : \epsilon \in E \} \text{ and } N_\epsilon = \sum_{i=1}^n I_{B_\epsilon}(X_i).$$

 We refer to the density of $\medidaAleatoria$ with respect to $\lambda$ by $\theta$ and to the posterior distribution induced by $\theta$ by $\posterior$. Furthermore, expectations with regard to $\posterior$ are denoted by $\mathbb{E}[\cdot | \Amostra]$ and expectations with regard to $f_0$ are denoted by $\ex[\cdot]$. The random variable $N_\epsilon = \sum_{i=1}^n I_{B_\epsilon}(X_i)$ and the realization $n_\epsilon$ are used, respectively, when analyzing the behavior of $\posterior$ as a random variable or a fixed realization. Realizations of $\posterior$ are important as they're probability measures over the set of densities with support on $\mathcal{X}$ that capture the updating of $\prior$ as a belief measure under the observation of $\Amostra$.

We restrict attention to specific choices of parameters for Pólya Trees. This ensures both the existence of $\Pi$ and the smoothness of random measures $\medidaAleatoria \sim \Pi$. Our results concern Pólya Trees $\ArvoreDePolya$ that satisfy:

\begin{assumption}[Polya tree hyperparameters]
 \label{ass:polya_hyper} \
  \begin{enumerate}
        \item If
        $l(\epsilon) = l(\epsilon^*)$, then $\alpha_\epsilon = 
        \alpha_{\epsilon^*}$. For simplicity, we define
        $a_{\len} := \alpha_{\epsilon}$.
        \item 
        \begin{equation}
    \sum_{l=1}^\infty \frac{l}{a_l} < \infty\tag{prior is smooth}
\end{equation}
        \item The partition structure  $\mathcal{B}$ satisfies $\lambda(\BE) = 2^{-\len}$.
        %\item $\mathcal{B}$ is the set of dyadic partitions of $[0,1]$, under which $B_\epsilon = \lambda^{-1}\left(\left[\sum \frac{\epsilon_i}{2^i}, \sum \frac{\epsilon_i}{2^i} + \frac{1}{2^i}\right)\right)$.
    \end{enumerate}
\end{assumption}

Under these conditions, $\medidaAleatoria$ is absolutely continuous with respect to $\lambda$ \citep{Ghosal_vanderVaart_2017}. If $\mathcal{X}=[0,1]$ then $\particao$ can be uniquely specified as the dyadic intervals of unity:

$$\mathcal{D} = \Biggl\{\left[\sum_{i=1}^n \frac{\epsilon_i}{2^i}, \sum_{i=1}^n \frac{\epsilon_i}{2^i} + \frac{1}{2^\len}\right], \epsilon \in E \Biggl\}.$$

For a general $\mathcal{X}$ the are many possible choices of $\particao$ and our findings are articulated for any one that satisfies $\lambda(B_\eps) = 2^{-\len}$.  Those are defined as \textbf{canonical} partitions of $\mathcal{X}$ with respect to the Lebesgue measure.

The restriction on $a_l$ in Assumption \ref{ass:polya_hyper} is required to ensure that $\theta$ samples from smooth densities. For instance, $\sum_{\eps \in E}\frac{1}{a_l} < \infty$ is a necessary and sufficient condition for $\medidaAleatoria$ to be absolutely continuous with respect to the Lebesgue measure \citep{Ghosal_vanderVaart_2017} and, hence, for $\theta$ to exist. By further requiring that $\sum_{\eps \in E} \frac{l}{a_l} < \infty$, $\theta$ is bounded away from $0$ and from $\infty$, as provided in the following novel result. This result ensures that the integral $\int _\mathcal{X}|\log\theta (t)|dt$ is almost surely bounded.

\begin{lemma}
 \label{lem:pt-ba0i}
 Under Assumption \ref{ass:polya_hyper},
 $\theta$ is almost surely bounded away 
 from $0$ and $\infty$.
\end{lemma}

\subsection{Regularity conditions on \texorpdfstring{$f_0$}{f0}}

The Theorems outlined in the previous section require some conditions on $f_0$, the data-generating distribution. First, we require that $f_0$ is in the KS-support of $\prior$, that is, there exists a positive probability that $\theta$ is drawn in a KS-neighborhood of $f_0$. This condition that is commonly required. 
Formally, $f_0$ is in the KS-support of $\theta$ if: 
\begin{assumption}[KS-support]
 \label{ass:ks_sup}
 There exists $\epsilon > 0$, such that
 $\prior(\theta: K(f_0, \theta) \leq \epsilon) > 0$.
\end{assumption}

In order to obtain posterior consistency , Assumption \ref{ass:ks_sup} is sufficient. For contraction rates, we require $f_0$ to be in a class of Hölder continuous densities on $[0,1]$ bounded away from $0$:
\begin{assumption}[Hölder class and bounded away from $0$]
 \label{ass:f0_holder} For some $m,K \in \mathbb{R^+}$ and  $\alpha \in[0,1]$:
 \begin{itemize}
  \item $f(t) > m \text{ for }t\in[0,1]$;
  \item $|f(x)-f(y)| \leq K|x-y|^\alpha  \text{ for all }x,y \in [0,1]$.
  \end{itemize}
\end{assumption}
Assumption \ref{ass:f0_holder} is standard for minimax results on density estimation. It is also appropriate for differential entropy estimation \citep{HallMorton1993, castillo_bernsteinvon_2015}. 

The consistency of our differential entropy estimator requires a weaker condition than Assumption \ref{ass:f0_holder}. Let
$y_{f_0}$ be
\begin{align*}
 y_{f_0}(\eps) &= \frac{P_{f_0}(B_{\eps1 ... \eps_k})}{P_{f_0}(B_{\eps1 ... \eps_{k-1}})}.
\end{align*}

One might expect that for $\posterior$ to concentrate at appropriate speed around $f_0$, then $Y_{\eps}$ should converge uniformly to $y_{f_0}(\eps)$.
Indeed, this is the case and for such a convergence to occur, Assumption \ref{ass:f0_holder} %guarantees that $y_{f_0}(\eps)$ is bounded above from $0$, bounded below from $1$, and converges to $0.5$ as $l(\eps) \rightarrow \infty$.
requires a very regular structure on $f_0$. The consistency of our differential entropy estimator, however, requires only the weaker condition that $y_{f_0}(\eps)$ is bounded from $1$. Whenever there is no ambiguity, we refer to $y_{f_0}(\epsilon)$ by $y_{\epsilon}$.
\begin{assumption}[Bounded $y^*_{f_0}$]
 \label{ass:f0_bounded}
 $y^*_{f_0} := \sup_{\eps} y_{f_0}(\eps) < 1$.
\end{assumption}

Besides those conditions, we also require $f_0$ to have finite differential entropy.

\begin{assumption}[Finite differential entropy]
 \label{ass:dif_entropy}
 \begin{align*}
 |H(f_0)| &:= 
 \left|\int_{\mathcal{X}} f_0(t) \log f_0(t) dt\right| < \infty
 \end{align*}
\end{assumption}

This condition is connected with the notion of Kullback-Leibler support. For Pólya Trees that satisfy Assumption \ref{ass:polya_hyper}, every $f_0$ that satisfies Assumption \ref{ass:dif_entropy} also satisfies Assumption \ref{ass:ks_sup}:

\begin{prop}[\citet{Ghosal_vanderVaart_2017}]
 Under Assumption \ref{ass:polya_hyper},
 Assumption \ref{ass:dif_entropy} implies
 Assumption \ref{ass:ks_sup}.
\end{prop}

%Proposition 1 states that an absolutely continuous Pólya Tree assigns mass to the neighborhood of any finite differential entropy density, but we require a stronger result. Theorem 1 is the result of being able to exchange the expectation and the product of the following identity, so we look for an application of the Dominated Convergence Theorem. %Specifically, we consider the case where $\left|\log \theta(t)\right|$ is bounded.

\section{Main results}\label{sec:main-results}

Next, the statements of the main Theorems in the paper are formally presented and discussed. Theorem 1 states that, under mild regularity conditions, Pólya Tree priors are strongly consistent in Kullback-Leibler neighborhoods.

\begin{theorem}[Posterior consistency]
 \label{teo:weak-consistency}
 Under Assumptions \ref{ass:polya_hyper} and \ref{ass:ks_sup},
 $\e[K(f_0, \theta)|\Amostra] \rightarrow 0$ almost surely and in in $L_1$. If additionaly
 Assumption \ref{ass:f0_holder} holds,
 $a_l = 2^{2\alpha l}$, and
 $\mathcal{X} = [0,1]$, then
 \begin{align*}
  \ex\left[\e[K(f_0, \theta)|\Amostra]\right] 
  = O\left(n^{-\frac{2\alpha}{1+2\alpha}}\right).
 \end{align*}
\end{theorem}

Previous works proved strong consistency in other semimetrics. \citet{schwartz_bayes_1965} establishes general conditions for consistency in the weak topology, which were  applied to Pólya Trees by \citet{lavine_aspects_1992}. Strong consistency with respect to the Hellinger and Total Variation metric was obtained by \citet{barron_consistency_1999} and further refined by \citet{walker_new_bayesian}. Those works obtained general bayesian nonparametric consistency theorems that were applied to specific Pólya Tree priors as examples. Consistency of Pólya Trees with respect to the supremum norm alongside convergence rates were obtained more recently by \citet{consistency_sup_norm}. All such results were stated in terms of one dimensional density estimation.

Theorem \ref{teo:weak-consistency} uses similar assumptions as the above papers. Assumption \ref{ass:ks_sup} is similar to the ones commonly used for obtaining strong consistency \citep{schwartz_bayes_1965,barron_consistency_1999}. 
Furthermore, the slowest known growth for Kullback-Leibler consistency was $a_l = l^{3+\delta}, \delta > 0$ \citep{walker_new_bayesian}. Assumption \ref{ass:polya_hyper} allows slower rates, such as 
$a_l = l^{2+\delta}$. Also, Assumption \ref{ass:f0_holder} has been used for obtaining the posterior contraction rate \citep{consistency_sup_norm}. 

Theorem \ref{teo:weak-consistency} also yields stronger conclusions than usual. To the best of our knowledge,
the convergence $\e[K(f_0, \theta)|\Amostra]$ had not yet been studied. This convergence guarantees robust consistency concerning Kullback-Leibler divergence, Total Variation and Hellinger distances. Moreover, the speed of convergence of $\e[K(f_0, \theta)|\Amostra]$ is a contraction rate for $\posterior$.

\begin{corollary}
 \label{cor:contract}
 Under Assumptions \ref{ass:polya_hyper} and \ref{ass:ks_sup}, 
 for $d \in \{d_{KL}, d_{H}, d_{TV}\}$,
 \begin{itemize}
   \item $\posterior$ is strongly consistent in $d$,
   \item For $\hat{\theta}(t) = \e\left[\theta(t)| \Amostra\right]$, $d(\hat{\theta}, f_0) \rightarrow 0$ almost surely and in $L_1$, and
   \item If Assumption \ref{ass:f0_holder} holds and
   $a_l = 2^{2\alpha l}$, then $\posterior$ 
   contracts at the minimax rate.
 \end{itemize}
\end{corollary}

To the best of our knowledge, Corollary \ref{cor:contract} establishes the first contraction rate for Pólya Trees with respect to the KL divergence. Furthermore, it establishes that, if $a_l = 2^{2\alpha l}$, then the posterior contracts at minimax rate.
Hence, $a_l = 2^{2\alpha l}$ obtains minimax contraction rates with respect to both the $L_1$ and supremum norms
\citep{castillo_bernsteinvon_2015}.

The convergence of $\e[K(f_0,\theta)|\Amostra]$ in Theorem \ref{thm:entropy} also implies the consistency of the \ref{eq:estimator-of-differential-entropy-alt}.

\begin{theorem}[Consistency of differential entropy estimator]
 \label{thm:entropy}
 Let $\hat{H}(\Amostra)$ be the \ref{eq:estimator-of-differential-entropy-alt}.
 Under Assumptions \ref{ass:polya_hyper}, \ref{ass:f0_bounded},  and \ref{ass:dif_entropy},
 if $a_l = 2^{\beta l}$ for $\beta >2$, then 
 $|\hat{H}(\Amostra)-H(f_0)| \rightarrow 0$ 
 in $\Amostra$ probability.
\end{theorem}

To the best of our knowledge, this is the first explicit estimator of differential entropy based on a full Bayesian Nonparametric posterior. Earlier Bayesian differential entropy methods considered finite alphabets or truncated priors \citep{nsb, castillo_bernsteinvon_2015}.

Theorems \ref{teo:weak-consistency} and \ref{thm:entropy} rely on a novel representation theorem for $K(f_0,\theta)$.

\begin{theorem}[Representation Theorem]
 \label{thm:representation}
 Under Assumptions \ref{ass:polya_hyper} and \ref{ass:ks_sup}, then
 \begin{equation}\label{eq:series-kl}
    K(f_0, \theta) = -H(f_0)-\sum_{\epsilon \in E} \left(P_X(B_{\epsilon0}) \log 2 Y_{\epsilon0} + P_X(B_{\epsilon1}) \log 2(1-Y_{\epsilon0})\right)\tag{series representation of KL divergence}
 \end{equation}
 \textit{Furthermore $\e\left[K(f_0,\theta)|X_1, \dots, X_n\right]$ is a supermartingale with respect to the filtration $\mathcal{F} = \{\sigma(X_1, \dots, X_n)\}_{n \in \mathbb{N}}$}.
\end{theorem}

Theorem \ref{thm:representation} expresses $K(f,\theta)$ as a sum involving the more tractable random variables $Y_{\epsilon 0}$. This insight plays a central role in the proofs of Theorems \ref{teo:weak-consistency} and \ref{thm:entropy}. Theorem \ref{thm:representation} also establishes that $K(f_0,\theta)$ is a supermartingale, that is, on average $\theta$ monotonically approaches $f_0$. Finally, Theorem \ref{thm:representation} also motivates the following differential entropy estimator:
\begin{align*}\label{eq:estimator-of-differential-entropy-alt}
 \estimador
 &= -\sum_{\epsilon \in E} \frac{\nne}{n} \mathbb{E}\left[\log(2Y_\epsilon)|\Amostra\right] = - \mathbb{E}\left[\sum_{i=1}^n \frac{\log \theta(X_i)}{n}|\Amostra\right]
 = \\
 &= -\frac{1}{n}\sum_{\epsilon \in E} \left(\nne\log(2) + \nze\psi(N_{\epsilon0}+a_l)+\nnu\psi(\nnu+a_l)-\nne\psi(\nne+2a_l)\right)\tag{differential entropy estimator}
\end{align*}

where $\psi$ denotes the digamma function.

Next, we discuss these results in detail and provide outlines for the proofs. Full proofs are provided in the appendix.

\subsection{Representation theorem}
\label{sec:KL-inc}

Theorem \ref{thm:representation} establishes that
$\e\left[K(f_0,\theta)|X_1, \dots, X_n\right]$ is a supermartingale with respect to the data filtration. Hence, on average, the separation between $\theta$ and $f_0$ decreases with the sample size. While the conclusion is expected, our proof leverages several unique properties of Pólya Tree priors. 

Theorem \ref{thm:representation} also shows that $K(f_0,\theta)+H(f_0)$ can be understood as a countable sum over the levels of the Pólya tree:
\begin{align*}
 K(f_0, \theta) + H(f_0) &= 
 -\sum_{\epsilon \in E} \left(P_X(B_{\epsilon0}) \log 2 Y_{\epsilon0} + P_X(B_{\epsilon1}) \log 2(1-Y_{\epsilon0})\right) \\
 &= -\sum_{j=1}^{\infty}\sum_{\epsilon \in E^*_j} \left(P_X(B_{\epsilon0}) \log 2 Y_{\epsilon0} + P_X(B_{\epsilon1}) \log 2(1-Y_{\epsilon0})\right).
\end{align*}
One of the advantages of this representation is its mathematical tractability. Note that $(Y_{\epsilon 0})_{\epsilon \in E}$ are jointly independent and appear in a single summand. This representation also indicates how well a truncated Pólya tree approximates the full tree. One might expect that the countable sum in Theorem \ref{thm:representation} can be truncated at a given level of the Polya tree and:
\begin{align}
 \label{eq:kl_finite}
 K(f_0, \theta) + H(f_0) &\approx 
 -\sum_{j=1}^{K}\sum_{\epsilon \in E^*_j} \left(P_X(B_{\epsilon0}) \log 2 Y_{\epsilon0} + P_X(B_{\epsilon1}) \log 2(1-Y_{\epsilon0})\right).
\end{align}
Indeed, Theorem \ref{teo:prob-cons} shows that the right side of eq. \ref{eq:kl_finite} converges a.s. to the left side. Furthermore, the right side of eq. \ref{eq:kl_finite} has been studied by \citet{Watson04052017}. It is the value of $K(f_0,\theta_K)+H(f_0)$, where $\theta_K$ is a Pólya tree truncated at length $K$.

Theorem \ref{thm:representation} plays a central role in the proof of Theorem \ref{teo:weak-consistency}. Using \citet{Ghosal_vanderVaart_2017}[Lemma B.10],
\begin{align}
 \label{eq:vaart}
 -H(f_0) &= \sum_{\epsilon \in E} \left(\fbo  \log\left(2\frac{\fbo}{\fbe}\right) + \fbu \log\left(2\frac{\fbu}{\fbe}\right)\right)
\end{align}
Applying eq. \ref{eq:vaart} to Theorem \ref{thm:representation}, one obtains
\begin{align*}
 K(f_0,\theta) 
 &= \sum_{\epsilon \in E} \left(\fbo  \log\left(\frac{\fbo}{\fbe} Y_{\eps 0}^{-1}\right) + \fbu \log\left(\frac{\fbu}{\fbe}Y_{\eps 0}^{-1}\right)\right)
\end{align*}
Hence, Theorem \ref{teo:weak-consistency} obtains $K(f,\theta) \rightarrow 0$ by showing that $Y_{\eps0} \rightarrow \frac{\fbo}{\fbe}$ uniformly over $\epsilon$.

% Both results shall be referenced throughout the paper and define the main restrictions we must impose on the prior $\prior$ and the true data generating process $f_0$.

The central argument in the proof of Theorem \ref{thm:representation} is one of ``expectation and limit commutes''. Let $X \sim f_0$ and consider its binary digits $0.\epsilon_1(X)\epsilon_2(X) \dots$. For a fixed realization of $\{Y_\epsilon\}_{\epsilon \in E}$,
\begin{align}
 \label{eq:classical-result-prior}
 K(f_0,\theta) + H(f_0)
 = \int f_0(t) \log \theta(t) dt
 &= \mathbb{E}_{X \sim f_0}\left[\log\theta(X)\right] \notag \\
 &= \mathbb{E}_{X \sim f_0}\left[\log\left(\prod_{j=1}^\infty 2 Y_{\epsilon_1(X) \dots \epsilon_l(X)}\right)\right] \notag \\  
 &= \mathbb{E}_{X \sim f_0}\left[\sum_{j=1}^\infty\log\left( 2 Y_{\epsilon_1(X) \dots \epsilon_l(X)}\right)\right].
\end{align}
The main argument in the proof of Theorem \ref{thm:representation} relies on showing that the expectation and the sum in eq. \ref{eq:classical-result-prior} commute and
\begin{align*}
 K(f_0,\theta) + H(f_0) 
 &= \sum_{j=1}^\infty \mathbb{E}_{X \sim f_0}\left[\log\left( 2 Y_{\epsilon_1(X) \dots \epsilon_l(X)}\right)\right].
\end{align*}
This passage is justified through the Dominated Convergence Theorem. Assumptions \ref{ass:polya_hyper} and \ref{ass:ks_sup} ensure that $\sup|\log(\theta(t))| < \infty$, which is an upper bound for $|\int f_0(t) \log \theta(t)dt|$.

\subsection{Consistency of Pólya Trees}

Theorem \ref{teo:prob-cons} establishes conditions that ensure $\e[K(f_0, \theta)|\Amostra] \stackrel{a.s.}{\rightarrow} 0$, thus yielding strong posterior consistency, since
\begin{align*}
 \posterior\left(K(f_0, \theta) > \epsilon \right)
 &\leq \frac{\e_{\posterior}[K(f_0, \theta)]}{\epsilon} 
 \stackrel{a.s.}{\rightarrow} 0.
\end{align*}
Furthermore, Theorem \ref{teo:prob-cons} also obtains conditions that ensure $\ex\left[\e_{\posterior}[K(f_0, \theta)]\right] = O\left(n^{-\frac{2\alpha}{1+2\alpha}}\right)$. This fact yields a posterior contraction rate, since
\begin{align*}
 P_{\Amostra}\left(\posterior\left(K(f_0, \theta) > M_n\epsilon_n \right) > \delta \right)
 &\leq P_{\Amostra}\left(
 \frac{\e_{\posterior}[K(f_0, \theta)]}{M_n\epsilon_n} > \delta \right) \\
 &\leq \frac{\ex\left[\e_{\posterior}[K(f_0, \theta)]\right]}{\delta M_n \epsilon_n}.
\end{align*}
As $O\left(n^{-\frac{2\alpha}{1+2\alpha}}\right)$ is the minimax rate of convergence for Hölder classes, this upper bound cannot be further improved. Furthermore, the above posterior contraction rate is obtained when $\al = 2^{2\alpha l(\eps)}$, which also attains the optimal posterior contraction rate according to the supremum norm \citep{castillo_bernsteinvon_2015}.

Polya Trees can obtain the minimax rate since Assumption \ref{ass:polya_hyper} places a high prior probability on adequately smooth densities. Indeed, Assumption \ref{ass:polya_hyper} is enough to ensure a soft truncation of rough elements of $\theta$. If $l(\epsilon)$ is sufficiently large, then the number of samples that fall in $B_{\epsilon}$ is at most $1$. Hence, since $\al \rightarrow \infty$, at high levels of the tree, the $Y_{\eps}$ are concentrated around $0.5$. The prior $\al = 2^{2\alpha \len}$ truncates the high levels of the tree at optimal speed.

The proof of Theorem \ref{teo:prob-cons} relies strongly on Theorem \ref{thm:representation}. The latter obtains that $K_n := \mathbb{E}\left[K(f_0, \theta)|\Amostra\right]$ is a supermartingale. Hence,
since $K(f_0, \theta) \geq 0$, it is sufficient to show that $\mathbb{E}_0[K_n] \rightarrow 0$ to establish that $K_n$ converges a.s. to $0$. Furthermore, Theorem \ref{thm:representation} yields an explicit representation of $K_n$ as the sum of independent terms that 
involve $Y_{\epsilon 0}$. This representation allows a bound over $\mathbb{E}_0[K_n]$.

\subsection{Differential entropy estimation}\label{sec:properties-estimator}

Theorem \ref{thm:entropy} obtains the consistency of the \ref{eq:estimator-of-differential-entropy-alt} requiring that $f_0$ satisfies Assumptions \ref{ass:f0_bounded} and \ref{ass:dif_entropy}.
Assumption \ref{ass:dif_entropy} requires solely that the parameter of interest, $H(f_0)$, is finite. Besides this Assumption, a common requirement for differential entropy estimators \citep{HallMorton1993} is that $f_0$ is bounded above $0$ and below $\infty$:
\begin{assumption}
 \label{ass:bounded_f0}
 There exists $M_1 > 0$ and $M_2 < \infty$ such that,
 for every $x \in \mathcal{X}$,
 $M_1 < f_0(x) < M_2$.
\end{assumption}
Assumption \ref{ass:bounded_f0} is stronger than Assumption \ref{ass:dif_entropy}. Indeed, under Assumption \ref{ass:bounded_f0},
$M_1 \lambda(B_\eps) < f_0(x) < M_2 \lambda(B_\eps)$. Hence,
\begin{align*}
 y_{f_0}(\eps_1\ldots\eps_k) 
 &= \frac{P_X(B_{\eps_1,\ldots,\eps_k})}{P_X(B_{\eps_1,\ldots,\eps_{k-1}})} \\
 &= \frac{P_X(B_{\eps_1,\ldots,\eps_k})}{P_X(B_{\eps_1,\ldots,\eps_k}) + P_X(B_{\eps_1,\ldots,(1-\eps_k)})} \\
 &\leq \frac{M_2 \lambda(B_{\eps_1,\ldots,\eps_k})}{M_2 \lambda(B_{\eps_1,\ldots,\eps_k}) + M_1 \lambda(B_{\eps_1,\ldots,(1-\eps_k)})} = \frac{M_2}{M_2+M_1} < 1.
\end{align*}

%First, Assumption \ref{ass:f0_bounded} is a simple requirement that may hold simply for any dimension and support, such as balls or continuous surfaces. 

%To have $f_0$ non vanishing is a common consideration because estimators mimic the arguably ideal estimator $T_n(f_0) = n^{-1}\sum \log f_0(X_i)$, not available in practice. A possible strategy for approximating $T_n(f_0)$ is to replace $f_0$ with an estimate $\hat{f}_n$, thereby obtaining a plug-in estimator $T_n(\hat{f_n})$. What harms this strategy is the observation of sample points in low density regions, under which $\hat{f}_n$ would typically be informed by few observations. Unlikely observations disturb $T_n(\hat{f_n})$ with outliers that may compromise convergence to $H(f_0)$. Regions in which $f_0$ grows too quickly may also negatively impact $T_n(\hat{f_n})$ for similar reasons.
Assumption \ref{ass:f0_bounded} allows $f_0$ to converge towards $\infty$, provided that it does not compromise the mass distribution across $\BO$ and $\BU$. That is, Assumption \ref{ass:f0_bounded} can be interpreted as a bound on the rate at which $f_0$ diverges to $\infty$.

\begin{example}
 If $f_0(x)$ follows the Beta(0.5,0.5) distribution, that is, $f_0(x) \propto x^{-0.5}(1-x)^{-0.5}$, for $x \in [0,1]$, then $f_0$ satisfies Assumption \ref{ass:f0_bounded} and does not satisfy Assumption \ref{ass:bounded_f0}.
\end{example}

\begin{example}
 If $f_0(x) \propto (e^{\sqrt{x}}-1)^{-1}$, for $x \in [0,1]$, then $f_0$ does not satisfy Assumption \ref{ass:f0_bounded}.
\end{example}

The proof of Theorem \ref{thm:entropy} relies on Theorems \ref{teo:prob-cons} and \ref{thm:representation}. The latter allows the following characterization of $H(f_0)$:
\begin{align*}
  H(f_0) &=
  \e_{\theta}[K(f_0, \theta)|\Amostra] - \sum_{\eps \in E}\fbe\e[\log(2Y_{\eps})|\Amostra]
\end{align*}
Furthermore, it follows from Theorem \ref{teo:prob-cons} that
$\e_{\theta}[K(f_0, \theta)|\Amostra]$ converges a.s. to 0. Hence,
\begin{align*}
 H^* := -\sum_{\eps \in E}\fbe\e[\log(2Y_{\eps})|\Amostra]\stackrel{a.s.}{\longrightarrow} H(f_0).
\end{align*}
Hence, it is sufficient to show that $|\estimador-H^*| \stackrel{a.s.}{\rightarrow} 0$ to complete the proof of Theorem \ref{thm:entropy}. By letting $\hat{P}_n(B_{\eps}) = \frac{\nne}{n}$, 
\begin{align*}
 |\estimador - H^*| &=
 \bigg|- \sum_{\eps \in E}\hat{P}_n\e[\log(2Y_{\eps})|\Amostra] - \sum_{\eps \in E}\fbe\e[\log(2Y_{\eps})|\Amostra]\bigg| \\
 &\leq \sup_{B_{\eps}}|\hat{P}_n(B_\eps)-\fbe| \cdot \bigg|\sum_{\eps \in E} \e[\log(2Y_{\eps})|\Amostra] \bigg|
\end{align*}
It follows from Glivenko-Cantelli that $\sup_{B_{\eps}}|\hat{P}_n(B_\eps)-\fbe| = O_P\left(n^{-0.5}\right)$.
Hence, it is enough to show that 
$\sum_{\eps}\e[\log(2Y_{\eps})|\Amostra] = o(n^{0.5})$, which follows from the choice of the prior $a_l = 2^{\beta l}$ with $\beta > 2$.

\section{Conclusion}\label{sec:conclusion}

This paper develops new analytical tools for studying Pólya Trees and related priors. Among our findings, we show that most Pólya Trees are almost surely bounded, and we introduce properties of $KL(f_0,\theta)$ that contribute to a more detailed understanding of their asymptotic behavior.

Our consistency results extend previous optimality findings for Pólya Tree priors with $\al = 2^{2\beta l(\eps)}$. Earlier work had established optimal contraction rates in the supremum norm for this class; here we show that the same sequence yields optimal contraction in both $L_1$ and Kullback–Leibler. We also demonstrate that this prior performs well for differential entropy estimation.

In the context of differential entropy estimation, the analysis establishes a link between two previously distinct literatures: the extensive work on nonparametric entropy estimation, including developments of the Kozachenko–Leonenko estimator, and the Bayesian nonparametric literature on density estimation. Although many Bayesian consistency results involve entropy functionals implicitly, explicit estimators of differential entropy under Bayesian nonparametric priors remain rare. Our results suggest that the regularity conditions commonly assumed for Bayesian density estimation can also yield effective, and in some cases optimal, entropy estimators.

Finally, several arguments developed here may extend directly to broader classes of priors, including Spike-and-Slab and Rubberly Pólya Trees. More generally, the proposed framework provides an analytically tractable approach to establishing consistency and contraction results for Pólya Tree–type priors. This approach can be viewed as a complement to general asymptotic theories based on Schwarz’s theorem, offering results that might be both stronger and more closely aligned with modern Bayesian convergence analyses.





% \end{thebibliography}

\appendix

\section{Proof of Lemma \ref{lem:pt-ba0i}}

The proof of Lemma \ref{lem:pt-ba0i} relies on auxiliary results that involve the moments of Beta random variables. These results are presented in Propositions \ref{prop:subgaussianity}, \ref{prop:beta_ineq}, and \ref{prop:beta_ineq_2}.

\begin{prop}\cite{subgaussiany}
 \label{prop:subgaussianity}
 Let $X \sim Beta(a, a)$ for $a \geq 1$. Then
 $$\mathbb{E}\left[\left(Y-\frac{1}{2}\right)^{2j}\right] = \frac{(2j)! (a)_j}{2^{2j}j! (2 a)_{2j}},
 \text{ where } (a)_j = \Gamma(a+j)/\Gamma(a).$$
\end{prop}

\begin{prop}
 \label{prop:beta_ineq}
 If $\theta \sim \ArvoreDePolya$, then
 \begin{align*}
  \e\left[\left(Y_\epsilon-\frac{1}{2}\right)^{2j}\right] \leq \frac{2 j^{j+1} e^{-j}}{(2\al + 1)^j}.
 \end{align*}
\end{prop}

\begin{proof}
 From Proposition \ref{prop:subgaussianity} we have
 
 $$\mathbb{E}_\theta\left[\left(Y_\epsilon-\frac{1}{2}\right)^{2j}\right] = \frac{(2j)! (\al)_j}{2^{2j}j! (2 \al)_{2j}}.$$
 
 We upper bound this number in two steps. First by the definition,
 
 $$\frac{(\al)_j}{(2 \al)_{2j}} \leq \frac{1}{2^j \prod_{k=1}^j (2\al + 2k-1)} \leq \frac{1}{2^j (2\al+1)^j}.$$
 
 Also, it follows from Stirling's inequality:

 $$\frac{(2j)!}{2^{2j}j!} \leq \frac{e (2j)^{2j+1/2} e^{-2j}}{\sqrt{2} (j)^{j+1/2} e^{-j}} \leq \frac{2^{2j+1} j^{2j+1} e^{-2j}}{j^{j} e^{-j}} = 2^{j+1} j^{j+1} e^{-j}.$$
 
 Therefore,
 
 $$\mathbb{E}_\theta\left[\left(Y_\epsilon-\frac{1}{2}\right)^{2j}\right] \leq \frac{2^{j+1} j^{j+1} e^{-j}}{2^j (2\al+1)^j} = \frac{2 j^{j+1} e^{-j}}{(2\al + 1)^j}.$$
\end{proof}

\begin{prop}
 \label{prop:beta_ineq_2}
 If $\theta \sim \ArvoreDePolya$, then there exists $C$ such that
 \begin{align*}
  \e\left[\max\limits_{\epsilon' \in E_{j-1}^*} \left(\left|\frac{1}{2}-Y_{\epsilon'0}\right|\right)\right] 
  \leq C\frac{j}{a_j}.
 \end{align*}
\end{prop}

\begin{proof}
 For every even valued $p$,
 \begin{align*}
  \e\left[\max\limits_{\epsilon' \in E_{j-1}^*} \left(\left|\frac{1}{2}-Y_{\epsilon'0}\right|\right)\right] 
  &\leq \e\left[\max\limits_{\epsilon' \in E_{j-1}^*} \left(\left|\frac{1}{2}-Y_{\epsilon'0}\right|\right)^p\right]^{1/p} & \text{Jensen's inequality} \\
  &= \e\left[\max\limits_{\epsilon' \in E_{j-1}^*} \left(\left|\frac{1}{2}-Y_{\epsilon'0}\right|^p\right)\right]^{1/p} \\
  &\leq \e\left[\sum_{\epsilon' \in E_{j-1}^*} \left(\left|\frac{1}{2}-Y_{\epsilon'0}\right|^p\right)\right]^{1/p} \\
  &= \left(2^{j-1} \e \left[\left|\frac{1}{2}-Y_{\epsilon'0}\right|^p \right]\right)^{\frac{1}{p}} \\
  &\leq \left(2^{j-1} \frac{2 p^{p+1} e^{-p}}{(2a_j + 1)^p}\right)^{\frac{1}{p}}
  & \text{Proposition \ref{prop:subgaussianity}}
 \end{align*}
 The proof is complete by choosing $p$ as 
 the first even value above $j$.
\end{proof}

\begin{proof}[Proof of Lemma \ref{lem:pt-ba0i}]
 Let $t = 0.\eps_1(t)\eps_2(t)\ldots$. For every $t$, $\theta(t) = \prod_{j=1}^{\infty}2Y_{\eps_1(t)\ldots\eps_j(t)}$. Therefore,
 for every $t$,
 \begin{align}
  \label{eq:pt-ba0i_1}
  \prod_{j=1}^\infty \min_{\epsilon \in \eje} 2Y_{\epsilon} \leq &\theta(t) \leq \prod_{j=1}^\infty \max_{\epsilon \in \eje} 2Y_{\epsilon}, \text{ and} \nonumber \\
  \prod_{j=1}^\infty \min_{\epsilon \in \eje} 2Y_{\epsilon} \leq \inf \theta(t) &< \sup \theta(t) \leq \prod_{j=1}^\infty 2\max_{\epsilon \in \eje} Y_{\epsilon}.
 \end{align}
 The right and left terms can be developed in the following way:
 \begin{align*}
  \max\limits_{\epsilon \in \eje} 2Y_{\epsilon} 
  &= \max\limits_{\epsilon' \in E_{j-1}^*} \{\max\{2 Y_{\epsilon'0}, 2Y_{\epsilon'1}\}\} \\
  &= \max\limits_{\epsilon' \in E_{j-1}^*} \left(1 + 2\left|\frac{1}{2}-Y_{\epsilon'0}\right|\right) \\
  &= 1+ 2\max\limits_{\epsilon' \in E_{j-1}^*} \left(\left|\frac{1}{2}-Y_{\epsilon'0}\right|\right) \text{, and } \\
  \min\limits_{\epsilon \in \eje} 2Y_{\epsilon} 
  &= 1- 2\max\limits_{\epsilon' \in E_{j-1}^*} \left(\left|\frac{1}{2}-Y_{\epsilon'0}\right|\right).
\end{align*}
 Hence, by applying the above inequalities in eq. \ref{eq:pt-ba0i_1},
 \begin{align}
  \label{eq:pt-ba0i_2}
  \prod_{j-1}^{\infty}\left(1- 2\max\limits_{\epsilon' \in E_{j-1}^*} \left(\left|\frac{1}{2}-Y_{\epsilon'0}\right|\right)\right) \leq \inf \theta(t) &< \sup \theta(t) \leq \prod_{j=1}^\infty \left(1+ 2\max\limits_{\epsilon' \in E_{j-1}^*} \left(\left|\frac{1}{2}-Y_{\epsilon'0}\right|\right)\right).
 \end{align}
 It remains to show that the left term in eq. \ref{eq:pt-ba0i_2} is positive and the right term is finite. Both statements are equivalent to
 \begin{align}
  \label{eq:pt-ba0i_3}
  \sum_{j=1}^\infty 2\max\limits_{\epsilon' \in E_{j-1}^*} \left|\frac{1}{2}-Y_{\epsilon'0}\right| < \infty
  \text{ a.s.}
 \end{align}
 Eq. \ref{eq:pt-ba0i_3} follows directly from observing that
 \begin{align*}
  \e\left[\sum_{j=1}^\infty2\max\limits_{\epsilon' \in E_{j-1}^*} \left(\left|\frac{1}{2}-Y_{\epsilon'0}\right|\right)\right]
  &= 2\sum_{j=1}^\infty\e\left[\max\limits_{\epsilon' \in E_{j-1}^*} \left(\left|\frac{1}{2}-Y_{\epsilon'0}\right|\right)\right]
  & \text{Monotone Convergence Theorem} \\
  &\leq  2C\sum_{j=1}^\infty \frac{j}{a_l}
  & \text{Proposition \ref{prop:beta_ineq_2}} \\
  & < \infty & \text{Assumption \ref{ass:polya_hyper}}
 \end{align*}

\end{proof}

\section{Proof of Theorem \ref{thm:representation}}

\begin{proof}[Proof of Theorem \ref{thm:representation}]

Let $X\sim f_0$ and consider $\epsilon_l(X)$ the $l$-th digit of its dyadic expansion. First we observe that 

$$\int f_0(t) \log \theta(t) dt = \mathbb{E}_{X \sim f_0}\left[\log\left(\prod_{l=1}^\infty 2 Y_{\epsilon_1(X) \dots \epsilon_l(X)}\right)\right] = \mathbb{E}_{X \sim f_0}\left[\sum_{l=1}^\infty \log\left(2 Y_{\epsilon_1(X) \dots \epsilon_l(X)}\right)\right].$$

However for any fixed $x \in [0,1]$ 
 
 $$\left|\sum_{l=1}^\infty \log\left(2 Y_{\epsilon_1(x) \dots \epsilon_l(x)}\right)\right| = |\log \theta (x)|$$

 and therefore

$$\left|\sum_{l=1}^\infty \log\left(2 Y_{\epsilon_1(X) \dots \epsilon_l(X)}\right)\right| < \sup_{x \in [0,1]} |\log \theta(x)|$$

which is $\prior$-almost surely bounded by Lemma 1. Thus, by the Dominated Convergence Theorem for almost all $\{Y_\epsilon\}_{\epsilon \in E}$,

$$\int f_0(t) \log \theta(t) dt = \mathbb{E}_{X \sim f_0}\left[\sum_{l=1}^\infty \log\left(2 Y_{\epsilon_1(X) \dots \epsilon_l(X)}\right)\right] = \sum_{l=1}^\infty \e_{X\sim f_0} \left[\log(2Y_{\eps_1(X) \dots \eps_l(X)})\right].$$

It follows that

\begin{align*}
    \int f_0(t) \log \theta(t) dt = \sum_{l=1}^\infty\mathbb{E}_{X \sim f_0}\left[\log\left(2 Y_{\epsilon_1(X) \dots \epsilon_l(X)}\right)\right] = \\ \sum_{l=1}^\infty \left(\sum_{\epsilon \in E_l^*} F_0(B_{\epsilon}) \log\left(2 Y_\epsilon\right)\right)= \sum_{\epsilon \in E} F_0(B_{\epsilon}) \log\left(2 Y_\epsilon\right).
\end{align*}

We conclude the proof of the first part of the Theorem noting that each $Y_\epsilon$ for $\len \geq  1$ is either $Y_{\epsilon'0}$ or $Y_{\epsilon'1} = 1-Y_{\epsilon'0}$ for a $\epsilon' \in E$. Thus:
$$\sum_{\epsilon \in E} F_0(B_{\epsilon}) \log\left(2 Y_\epsilon\right)  = \sum_{\epsilon \in E} \left(F_0(B_{\epsilon0}) \log\left(2 Y_{\epsilon0}\right) + F_0(B_{\epsilon1}) \log\left(2 (1-Y_{\epsilon0})\right)\right) $$

which is a sum of mutually independent random variables.

For the second part of the Theorem 2 let $K_n = \e\left[K(f_0, \theta)|\Amostra\right]$. We have just proved that
$$K_n = \int f(t) \log f(t)dt-\mathbb{E}_{\theta}\left[\sum_{\epsilon \in E} \left(P_X(B_{\epsilon0}) \log 2 Y_{\epsilon0} + P_X(B_{\epsilon1}) \log 2(1-Y_{\epsilon0})\right)|\Amostra\right]$$

and consequently $K_n \geq 0$ for all $n$. Also, as the terms of the infinite sum are always negative by the monotone convergence theorem we have
\begin{align*}
K_n & =  \int f(t) \log f(t)dt-\sum_{\epsilon \in E} \mathbb{E}_{Y}\left[\left(P_X(B_{\epsilon0}) \log 2 Y_{\epsilon0} + P_X(B_{\epsilon1}) \log 2(1- Y_{\epsilon0})\right)|\Amostra\right] = \\
& -H(f)  \\
& -\sum_{\epsilon} \fbo\left(\log(2) + \psi(\nze+\al) -\psi(\nne+2\al)\right) \\
& -\sum_{\epsilon} \fbu\left(\log(2) + \psi(\nnu+\al)-\psi(\nne+2\al)\right)
\end{align*}

Now we prove that $K_n$ is a supermartingale. First let
\begin{align*}
T_{n,\epsilon} = & \fbo\left(\log(2) + \psi(\nze+\al) -\psi(\nne+2\al)\right) + \\
& \fbu\left(\log(2) + \psi(\nnu+\al)-\psi(\nne+2\al)\right)
\end{align*}

thus $K_n  = -H(f_0,f_0) - \sum_{\eps} T_{n,\eps}$ almost surely for all $n$. Now defining for $n \geq 1$
\begin{align*}
T_{n-1, \epsilon}^\leftarrow :=&\fbo\left(\psi\left(\sum_{i=1}^{n-1} \mathbb{I}_{\BO}(X_i) + \al+1\right) - \psi\left(\sum_{i=1}^{n-1} \mathbb{I}_{\BE}(X_i) + 2\al+1\right)\right) &+\\
&\fbu\left(\psi\left(\sum_{i=1}^{n-1} \mathbb{I}_{\BU}(X_i) + \al\right) - \psi\left(\sum_{i=1}^{n-1} \mathbb{I}_{\BE}(X_i) + 2\al+1\right)\right) &= \\
& \fbo\left(\psi\left(\sum_{i=1}^{n-1} \mathbb{I}_{\BO}(X_i) + \al\right) + \frac{1}{\sum_{i=1}^{n-1} \mathbb{I}_{\BO}(X_i) + \al}\right)& +\\
& \fbo\left(-\psi\left(\sum_{i=1}^{n-1} \mathbb{I}_{\BE}(X_i) + 2\al\right)-\frac{1}{\sum_{i=1}^{n-1} \mathbb{I}_{\BE}(X_i) + 2\al}\right) & + \\ 
&  \fbu\psi\left(\sum_{i=1}^{n-1} \mathbb{I}_{\BU}(X_i) + \al\right) & + \\ 
&  \fbu\left(- \psi\left(\sum_{i=1}^{n-1} \mathbb{I}_{\BE}(X_i) + 2\al\right)-\frac{1}{\sum_{i=1}^{n-1} \mathbb{I}_{\BE}(X_i) + 2\al}\right) &= \\ 
& T_{n-1,\eps} + \frac{\fbo}{\sum\limits_{i=1}^{n-1} \mathbb{I}_{\BO}(X_i) + \al}-\frac{\fbe}{\sum\limits_{i=1}^{n-1} \mathbb{I}_{\BE}(X_i) + 2\al}
\end{align*}

and analogously

\begin{align*}
T_{n-1, \epsilon}^\rightarrow := &\fbo\left(\psi\left(\sum_{i=1}^{n-1} \mathbb{I}_{\BO}(X_i) + \al\right) - \psi\left(\sum_{i=1}^{n-1} \mathbb{I}_{\BE}(X_i) + 2\al+1\right)\right) &+\\
&\fbu\left(\psi\left(\sum_{i=1}^{n-1} \mathbb{I}_{\BU}(X_i) + \al+1\right) - \psi\left(\sum_{i=1}^{n-1} \mathbb{I}_{\BE}(X_i) + 2\al+1\right)\right) &= \\
& T_{n-1,\eps} + \frac{\fbu}{\sum\limits_{i=1}^{n-1} \mathbb{I}_{\BU}(X_i) + \al}-\frac{\fbe}{\sum\limits_{i=1}^{n-1} \mathbb{I}_{\BE}(X_i) + 2\al} &
\end{align*}

We prove that $K_n$ is a supermartingale by exploring the following recursive relation:

$$T_{n,\epsilon} = \mathbb{I}_{\BO}(X_n)T_{n-1,\epsilon}^\leftarrow+ \mathbb{I}_{\BU}(X_n)T_{n-1,\epsilon}^\rightarrow +   \mathbb{I}_{\BE^c}(X_n)T_{n-1,\epsilon}$$

In words, observing a new sample point $X_n$ may not change the value of $T_{n,\eps}$ depending on wether $X_n$ falls in $\BE$. If it does, the digamma function allow for a explicit computation of the difference. Averaging over $X_n$ and conditioning on $\mathbf{X}_{n-1}$ we have

$$E[T_{n, \epsilon}|X_1, \dots X_{n-1}] = P_X(\BO)T_{n-1,\epsilon}^\leftarrow + P_X(\BU)T_{n-1,\epsilon}^\rightarrow + (1-\fbe) T_{n-1,\epsilon} = $$

$$T_{n-1, \epsilon}+\fbo\left(T_{n-1,\epsilon}^\leftarrow-T_{n-1,\epsilon}\right) +\fbu\left(T_{n-1,\epsilon}^\rightarrow-T_{n-1,\epsilon}\right).$$

Therefore

\begin{align*}E[T_{n, \epsilon}|X_1, \dots X_{n-1}] = & T_{n-1,\epsilon} & +  \\ &  \frac{\fbu^2}{\sum\limits_{i=1}^{n-1} \mathbb{I}_{\BU}(X_i) + \al}+ \frac{\fbo^2}{\sum\limits_{i=1}^{n-1} \mathbb{I}_{\BO}(X_i) + \al} &-\\ & \frac{\fbe^2}{\sum\limits_{i=1}^{n-1} \mathbb{I}_{\BE}(X_i) + 2\al} &
\end{align*}

We conclude noting that

\begin{align*}
& \frac{\fbu^2}{\sum\limits_{i=1}^{n-1} \mathbb{I}_{\BU}(X_i) + \al}+\frac{\fbo^2}{\sum\limits_{i=1}^{n-1} \mathbb{I}_{\BO}(X_i) + \al}-\frac{\fbe^2}{\sum\limits_{i=1}^{n-1} \mathbb{I}_{\BE}(X_i) + 2\al} & = \\
& \frac{\fbu^2}{\sum\limits_{i=1}^{n-1} \mathbb{I}_{\BU}(X_i) + \al}+\frac{\fbo^2}{\sum\limits_{i=1}^{n-1} \mathbb{I}_{\BO}(X_i) + \al}-\frac{(\fbu+\fbo)^2}{\sum\limits_{i=1}^{n-1} \mathbb{I}_{\BE}(X_i) + 2\al} & = \\
&\fbu^2 \left(\frac{\sum\limits_{i=1}^{n-1} \mathbb{I}_{\BO}(X_i) + \al}{\left(\sum\limits_{i=1}^{n-1} \mathbb{I}_{\BU}(X_i) + \al\right)\left(\sum\limits_{i=1}^{n-1} \mathbb{I}_{\BE}(X_i) + 2\al\right)}\right) & + \\
&\fbo^2\left(\frac{\sum\limits_{i=1}^{n-1} \mathbb{I}_{\BU}(X_i) + \al}{\left(\sum\limits_{i=1}^{n-1} \mathbb{I}_{\BO}(X_i) + \al\right)\left(\sum\limits_{i=1}^{n-1} \mathbb{I}_{\BE}(X_i) + 2\al\right)}\right) & - \\
&2 \fbu \fbo \left(\frac{1}{\sum\limits_{i=1}^{n-1} \mathbb{I}_{\BE}(X_i) + 2\al}\right) &  = \\
&\frac{\left(\fbo\left(\sum\limits_{i=1}^{n-1} \mathbb{I}_{\BO}(X_i)+\al\right)-\fbu\left(\sum\limits_{i=1}^{n-1} \mathbb{I}_{\BU}(X_i)+\al\right)\right)^2}{\left(\sum\limits_{i=1}^{n-1} \mathbb{I}_{\BE}(X_i) + 2\al\right)\left(\sum\limits_{i=1}^{n-1} \mathbb{I}_{\BU}(X_i) + \al\right)\left(\sum\limits_{i=1}^{n-1} \mathbb{I}_{\BO}(X_i) + \al\right)}&
\end{align*}

is positive which ensures $\mathbb{E}[K_n| X_1, \dots, X_{n-1}] \leq K_{n-1}$ as stated.

\end{proof}

\section{Proof of Theorem \ref{teo:weak-consistency}}

We present a proof that employs three auxiliary results alongside Theorem 3. Proposition \ref{prop:negative-moment} provides estimates for the negative moment of a binomial random variable and Proposition \ref{prop:digamma-inequality} provides bounds for evaluations of the digamma function [see e.g. \citep{gradshteyn2014table} for details]. Finally, Proposition \ref{prop:upper-bound-term} applies Proposition \ref{prop:negative-moment}, obtaining useful estimates for related quantities of $\e_{\Amostra}[\log(\hat{Y}_\epsilon)]$.

Throughout this session, we write $y_{\eps0} = \fbo/\fbe$, $y_{\eps1} = \fbu/\fbe$ and $c_{\eps} = \fbo \log\left(2y_{\eps0}\right)+\fbu \log\left(2y_{\eps1}\right)$

\begin{prop}[\citep{upper-bound-1overx}]\label{prop:negative-moment}
    Let $X \sim Bin(n,p)$ be a binomial random variable. It holds

    \begin{equation}
        \mathbb{E}\left[\frac{1}{a+X}\right] \leq \frac{1}{(np + a - 1+p)}
    \end{equation}
    
\end{prop}

\begin{prop}[\citep{batir_inequalities_2008}]\label{prop:digamma-inequality}

    The digamma function $\psi$ satisfy for all $x > 0$:

    \begin{equation}
\log(x)-\frac{1}{x} < \psi(x) < \log(x)-\frac{1}{2x}
\end{equation}
\end{prop}

\begin{prop}\label{prop:upper-bound-term}

    Under the regular Pólya Tree model, consider 

    $$T_{n,\eps} = \fbo \e\left[\log\left(2Y_{\eps0}\right)|\Amostra\right] + \fbu \e\left[\log\left(2Y_{\eps1}\right)|\Amostra\right].$$

    Then

    \begin{equation}\label{eq:limit}
        \mathbb{E}[T_{n,\epsilon}] \rightarrow c_\epsilon
    \end{equation}

    \begin{equation}\label{eq:bound1}
        |\mathbb{E}[T_{n,\epsilon}]| \leq \al^{-1}+c_\epsilon
    \end{equation}
    
\end{prop}

\begin{prop}\label{prop:bound-holder}

    If $\mathcal{X}=[0,1]$ and $f_0$ satisfies Assumption \ref{ass:f0_holder} there are $c_1$ and $c_2$ such that

    \begin{equation}\label{eq:bound3}
        c_\epsilon - \mathbb{E}[T_{n,\epsilon}]  \leq \frac{2^{\len}\fbe}{nc_1+2^{\len}\al}  + \frac{\al c_2}{2^{(1+2\alpha)\len}(nc_1+2^{\len}\al)}  
    \end{equation}
    
\end{prop}

We present our proof of Propositions \ref{prop:upper-bound-term} and \ref{prop:bound-holder} before proving Theorem \ref{teo:prob-cons} as they are a slight distraction from the main argument.

\begin{proof}[Proof of Proposition \ref{prop:upper-bound-term}]

First we observe that

$$\e[Y_{\eps0}|\Amostra] = \frac{\nze+\al}{\nne +2\al}$$

$$\e[Y_{\eps0}^{-1}|\Amostra] = \frac{\nne+2\al-1}{\nze +\al-1}$$

and that $\nze | \nne \sim Bin(\nne, y_{\eps0})$.

We prove equation \eqref{eq:limit} in two parts. First, applying Jensen's inequality twice, we obtain:
\begin{align*}
\ex[T_{n,\eps}]  = & \ex\left[\fbo \e\left[\log\left(2Y_{\eps0}\right)|\Amostra\right] + \fbu \e\left[\log\left(2Y_{\eps1}\right)|\Amostra\right]\right]  & \leq \\ 
&\fbo \log\left(2\ex\left[\frac{\nze+\al}{N_{\eps}+2\al}\right]\right) + \fbu \log\left(2\ex\left[\frac{\nnu+\al}{N_{\eps}+2\al}\right]\right)&  \leq \\
&\fbo \log\left(2\ex\left[\frac{\nze+\al}{N_{\eps}+1}\right]\right) + \fbu \log\left(2\ex\left[\frac{\nnu+\al}{N_{\eps}+1}\right]\right)  & \leq \\
&\fbo \log\left(2y_{\eps0}+2\ex\left[\frac{\al}{N_{\eps}+1}\right]\right) + \fbu \log\left(2y_{\eps1}+2\ex\left[\frac{\al}{N_{\eps}+1}\right]\right) &
\end{align*}

Thus $\ex\left[\frac{\al}{N_{\eps}+1}\right] \rightarrow 0$ by proposition \ref{prop:negative-moment} and also $\lim \ex[T_{n,\eps}] \leq c_\eps$. Now, we observe that Jensen's inequality also gives a lower bound:

\begin{align*}
\ex[T_{n,\eps}] = \ex\left[\fbo \e\left[\log\left(2Y_{\eps0}\right)|\Amostra\right] + \fbu \e\left[\log\left(2Y_{\eps1}\right)|\Amostra\right]\right] = \\
-\ex\left[\fbo \e\left[\log\left((2Y_{\eps0})^{-1}\right)|\Amostra\right] - \fbu \e\left[\log\left((2Y_{\eps1})^{-1}\right)|\Amostra\right]\right] \geq \\
%%$$\ex\left[\fbo \e\left[\log\left((2Y_{\eps0})^{-1}\right)\right] + \fbu \e\left[\log\left((2Y_{\eps1})^{-1}\right)\right]\right] \leq $$
-\fbo \log\left(\frac{1}{2}\ex\left[\frac{\nne+2\al-1}{\nze+\al-1}\right]\right) - 
\fbu \log\left(\frac{1}{2}\ex\left[\frac{\nne+2\al-1}{N_{\eps1}+\al-1}\right]\right).
\end{align*}

We also obtain

$$\ex\left[\frac{\nne+2\al-1}{\nze+\al-1}\right] \leq \ex\left[\frac{\nne+2\al-1}{\nne y_{\eps0}+1}\right] \leq \frac{1}{y_{\eps0}} + \ex\left[\frac{2\al}{\nne y_{\eps0}+1}\right]$$

by applying Proposition \ref{prop:negative-moment}. Thus, as $-\log(x)$ is a decreasing function it follow that

\begin{align*}
& -\fbo \log\left(\frac{1}{2}\ex\left[\frac{\nne+2\al-1}{\nze+\al-1}\right]\right) - \fbu \log\left(\frac{1}{2}\ex\left[\frac{\nne+2\al-1}{N_{\eps1}+\al-1}\right]\right) & \geq \\
& -\fbo \log\left(\frac{1}{2y_{\eps0}} + \frac{1}{2}\ex\left[\frac{2\al}{\nne y_{\eps}}\right]\right) - \fbu \log\left(\frac{1}{2y_{\eps1}} + \frac{1}{2}\ex\left[2\frac{2\al}{\nne y_{\eps}}\right]\right) & \rightarrow \\&c_\eps &
\end{align*}

Thus we conclude that $\ex[T_{n,\eps}] \rightarrow c_{\eps}$. For equation \eqref{eq:bound1} we note that by Theorem \ref{thm:representation} $\ex[T_{n,\eps}]$ is increasing. Therefore we have for all $n$

$$\ex[T_{n,\eps}] \leq \lim_{n \rightarrow \infty} \ex[T_{n,\eps}] = c_{\eps}.$$

However, as $T_{n,\eps}$ is a submartingale and applying Proposition \ref{prop:digamma-inequality}

$$\ex[T_{n,\eps}] \geq  -\al^{-1}$$

and an application of the triangle inequality concludes the proof.

\end{proof}

\begin{proof}[Proof of Proposition 8]

First we note that

\begin{align*}
& c_\eps-\ex[T_{n,\eps}] =  \\ 
& \ex\left[\fbo \e\left[\log\left(y_{\eps0}(Y_{\eps0})^{-1}\right)|\Amostra\right] + \fbu \e\left[\log\left(y_{\eps1}(Y_{\eps1})^{-1}\right)|\Amostra\right]\right] = \\
& \fbo\log\left(y_{\eps0}\ex\left[\frac{\nne+2\al-1}{\nze+\al-1}\right]\right) +\fbu\log\left(y_{\eps1}\ex\left[\frac{\nne+2\al-1}{\nnu+\al-1}\right]\right) \leq  \\
& \fbo\log\left(y_{\eps0}\ex\left[\frac{\nne+2\al-1}{\nne y_{\eps0}+\al-2+y_{\eps0}}\right]\right) +\\
&\fbu\log\left(y_{\eps1}\e_{X}\left[\frac{\nne+2\al-1}{\nne y_{\eps1}+\al-2+y_{\eps1}}\right]\right)  \leq \\
& \fbo\left(\ex\left[\frac{(2y_{\eps0}-1)\al+2-2y_{\eps0}}{\nne y_{\eps0}+\al-2+y_{\eps0}}\right]\right) + \\
&\fbu\left(\ex\left[\frac{(2y_{\eps1}-1)\al+2-2y_{\eps1}}{\nne y_{\eps1}+\al-2+y_{\eps1}}\right]\right) \leq \\
& \fbo\left(\ex\left[\frac{(2y_{\eps0}-1)\al+2}{\nne y_{\eps0}+\al-2}\right]\right) +\fbu\left(\ex\left[\frac{(2y_{\eps1}-1)\al+2}{\nne y_{\eps1}+\al-2}\right]\right) \leq \\
& \fbo\left(\frac{\al(2y_{\eps0}-1)+2}{n\fbo+\al-3+\fbo}\right) + \\
&\fbu\left(\frac{\al(2y_{\eps1}-1)+2}{n\fbu+\al-3+\fbu}\right) \leq \\
& \frac{2\fbe}{\al+n\min\{\fbu, \fbo\}} + \\
& \al \left(\frac{\fbo(2y_{\eps0}-1)}{n\fbo+\al-3}+\frac{\fbu(2y_{\eps1}-1)}{n\fbu+\al-3}\right) \leq \\
& \frac{2\fbe}{\al+n\min\{\fbu, \fbo\}-3}  + \frac{\al\frac{(\fbu-\fbo)^2}{\fbe}}{n\min\{\fbu, \fbo\}+\al-3}
\end{align*}

Now we observe that by the mean value theorem for every $\eps$ there are $x_{\eps0}$ and $x_{\eps1}$ such that

\begin{align*}
& \frac{(\fbu-\fbo)^2}{\fbe} = \frac{((f(x_{\eps0})-f(x_{\eps1})) 2^{-(\len+1)})^2}{(f(x_{\eps0})+f(x_{\eps1})) 2^{-\len}} \leq \\ & \frac{|x_{\eps0}-x_{\eps1}|^{2\alpha} 2^{-2\len}}{4m 2^{-\len}} \leq \frac{2^{-2\alpha(\len)}2^{-\len}}{4m}
\end{align*}

Where the last inequality arises from the fact that for any two points inside $\BE$, one must have $|x-y| \leq 2^{-\len}$, directly following from the Holder condition. Also

$$\min\{\fbu, \fbo\} = 2^{-(\len+1)}\min\{f(x_{\eps0}), f(x_{\eps1})\} > m2^{-\len-1}$$

By combining the last three inequalities, we obtain

$$c_\eps-\ex[T_{n,\eps}] \leq \frac{2\fbe}{\al+n m 2^{-\len-1}-3}  + \frac{\al\frac{2^{-2\alpha(\len+1)}2^{-\len}}{2m} }{nm2^{-\len-1}+\al-3}$$

and the result follows directly.

\end{proof}

\begin{proof}[Proof of Theorem \ref{teo:prob-cons}]

For the first part of the Theorem, the Martingale Convergence Theorem guarantees that $K_n \rightarrow K_\infty$ almost surely for some finite $K_\infty$ with $\e[K_\infty] \leq \e[K_n]$ for all $n$. We conclude the proof noting that equation \eqref{eq:bound1} ensures the summations are bounded and by Dominated Convergence Theorem we have

$$\lim_n\e[K_n] = \lim_n\sum_{\eps} \left(c_\eps-\e[T_{n,\eps}]\right) = \sum_{\eps} \lim_n\left(c_\eps-\e[T_{n,\eps}]\right) = 0.$$

For the second part of the Theorem we control the terms $c_\eps-\e[T_{n,\eps}]$ uniformly on $\eps$ to obtain an explicit upper bound for $\e[K_n]$. By Proposition \ref{prop:bound-holder}:

$$\e[K_n] = \sum_{\eps} \left(c_\eps-\e[T_{n,\eps}]\right)\leq \sum_{\eps}\left(\frac{2\fbe 2^{\len}}{n m/2+(\al-3)2^{\len}}  + \frac{\al 2^{-2\alpha \len} C }{nm/2+(\al-3)2^{\len}}\right) \leq $$

$$\sum_{l=1}^\infty\left(\frac{2^{l+1}}{n m/2+(a_l-3)2^{l}}  + \frac{2^{l}a_l 2^{-2\alpha l} C }{nm/2+(a_l-3)2^{l}} \right)$$

without loss of generality lets consider $a_l = 3 + 2^{2l\alpha}$. Then:

$$\e[K_n]  \leq \sum_{l=1}^\infty\left(\frac{22^l}{n m/2+2^{l(2\alpha+1)}}  + \frac{2^{l(2\alpha+1)} 2^{-2\alpha l} C }{nm/2+2^{l(2\alpha+1)}} \right).$$

Let $k_n = \log_2(n)(1+2\alpha)^{-1}$. Then

$$\sum_{l=1}^\infty \frac{2^{l+1}}{n m/2+2^{l(2\alpha+1)}} =\sum_{l=1}^{k_n} \left(\frac{2^{l+1}}{n m/2+2^{l(2\alpha+1)}}\right) + \sum_{l=k_n+1}^{\infty} \left(\frac{2^{l+1}}{n m/2+2^{l(2\alpha+1)}}\right) \leq $$

$$\sum_{l=1}^{k_n} \left(\frac{2^{l+1}}{n m/2+n}\right) + \sum_{l=k_n+1}^{\infty} \left(\frac{2^{l+1}}{2^{l(2\alpha+1)}}\right) = O(n^{-\frac{2\alpha}{1+2\alpha}})$$

and

$$\sum_{l=1}^\infty \left(\frac{2^{l(2\alpha+1)} 2^{-2\alpha l} C }{nm/2+2^{l(2\alpha+1)}}\right) \leq \sum_{l=1}^{k_n} \left(\frac{2^{l(2\alpha+1)} 2^{-2\alpha l} C }{nm/2+2^{l(2\alpha+1)}}\right) + \sum_{l=k_n}^\infty \left(\frac{2^{l(2\alpha+1)} 2^{-2\alpha l} C }{nm/2+2^{l(2\alpha+1)}}\right) = O(n^{-\frac{2\alpha}{1+2\alpha}}).$$

The last two inequalities guarantees that $\e[K_n] = O(n^{-\frac{2\alpha}{2\alpha+1}})$ as stated.

\end{proof}

\section{Proof of Theorem 2}

\newcommand\seps{s(\eps)}
\newcommand\sepsl{s(\epsl)}
\newcommand\nnepss{N_{\epss}}
\newcommand\nnepsso{N_{\epss0}}
\newcommand\expec{\ex[\nnel\parcelaRl \nne\parcelaR]}
\newcommand\pe{P(\BE)}
\newcommand\po{P(\BO)}
\newcommand\pel{P(B_{\epsl})}
\newcommand\pes{P(B_{\epss})}
\newcommand\pesu{P(B_{\epss1})}
\newcommand\peso{P(B_{\epss0})}

\newcommand\bineesu{\frac{\pe}{\pesu}}
\newcommand\bineeso{\frac{\pe}{\peso}}
\newcommand\binelesu{\frac{\pel}{\pesu}}
\newcommand\bineleso{\frac{\pel}{\peso}}
\newcommand\parcela{\log(2) + \psi(\al + \nze) - \psi(2\al + \nne)}

Our proof employs two auxiliary results.

\begin{lemma}
 \label{lemma:thm_2}
    Let $\theta \sim \ArvoreDePolya$ and $\theta|\Amostra$ its posterior. If

    $$a_l = 2^{l(2+\delta)}, \delta > 0$$

    then

    \begin{equation}
        \frac{\left(\sum_{\epsilon \in E}\ex\left[\left|\e\left[\log\left(2Y_\eps\right)|\Amostra\right]\right|\right]\right)}{\sqrt{n}} \rightarrow 0
    \end{equation}

    in $\Amostra$ probability.
    
\end{lemma}

\begin{prop}\label{prop:upper-bound-summand}

 Let $\theta \sim \ArvoreDePolya$ and $\theta|\Amostra$ be its posterior. Then for all $k\geq 2$
    \begin{equation}|\e\left[\log\left(2Y_{\eps_1 \dots \eps_k}\right)|\Amostra\right]| \leq \frac{N_{\eps_1 \dots \eps_{k-1}}}{\al}
    \end{equation}

    and, under Assumption \ref{ass:f0_bounded}, there is a constant $C > 1$ that depends on $f$ such that
    \begin{equation}\ex\left[|\e\left[\log\left(2Y_\eps\right)|\Amostra\right]|\right] \leq C
    \end{equation}

\end{prop}

\begin{proof}[Proof of Proposition \ref{prop:upper-bound-summand}]

Let $\eps = (\eps_1 \dots \eps_{k-1})$. Without loss of generality, let us assume that $\eps_k = 0$. First, we observe that

$$\log\left(\frac{2\nze+2\al-1}{\nne+2\al-1}\right) \leq \e\left[\log\left(2Y_{\eps0}\right)|\Amostra\right] \leq \log\left(\frac{2\nze+2\al}{\nne+2\al}\right)$$

which leads to 

$$\frac{2\nze-\nne}{\nze+2\al-1}\leq \e\left[\log\left(2Y_{\eps0}\right)|\Amostra\right] \leq \frac{2\nze-\nne}{\nne+2\al}$$

Therefore

$$|\e\left[\log\left(2Y_{\eps0}\right)|\Amostra\right]| \leq \frac{|\nze-\nnu|}{\nze
+2\al-1} \leq \frac{\nne}{\nze
+2\al-1} \leq \frac{\nne}{\al}$$

as stated. Also

$$|\e\left[\log\left(2Y_{\eps0}\right)|\Amostra\right]|\leq \frac{\nne}{\nze
+2\al-1} \leq \frac{\nne}{\nze+1}$$

and then by Proposition \ref{prop:negative-moment} we have

$$\e_X\left[\frac{\nne}{\nze+1}\right] \leq \ex\left[\frac{\nne}{\nne y_{\eps0}+y_{\eps0}}\right] < \frac{1}{\inf_\eps y_{\eps}}.$$

By Assumption \ref{ass:f0_bounded} the upper bound is finite and the result follows.

\end{proof}

\begin{proof}[Proof of Theorem \ref{thm:entropy}]

It holds

\begin{align*}|H(f_0) - \estimador| & \leq \\
\left|\estimador-\e\left[\int f_0(t) \log\theta(t)dt |\Amostra\right]\right| + \left|\e\left[\int f_0(t) \log\theta(t)dt |\Amostra\right]-H(f_0)\right| & = \\
\left|\estimador-\e\left[\int f_0(t) \log\theta(t)dt |\Amostra\right]\right| + \e\left[K(f_0; \theta)|\Amostra\right] &
\end{align*}

By Theorem 1 second term converges to $0$ almost surely. We conclude the proof showing that $\left|\hat{H}(\Amostra)-\e\left[\int f_0(t) \log\theta(t)dt |\Amostra\right]\right| \rightarrow 0$ in probability.

It holds

$$\left|\hat{H}(\Amostra)-\e\left[\int f_0(t) \log\theta(t)dt |\Amostra\right]\right| = $$

$$\left|\sum \frac{\nne}{n} \e\left[\log\left(2Y_\eps\right)|\Amostra\right]-\sum P_X(\BE)\e\left[\log\left(2Y_\eps\right)|\Amostra\right]\right| =$$

$$\left|\sum \left(\frac{\nne}{n}-P_X(\BE) \right)\e\left[\log\left(2Y_\eps\right)|\Amostra\right]\right| \leq \sum_{\epsilon \in E} \left|\frac{\nne}{n}-P_X(\BE) \right|\left|\e\left[\log\left(2Y_\eps\right)|\Amostra\right]\right| \leq $$
$$\left(\sup_{\epsilon \in E} \left|\frac{\nne}{n}-P_X(\BE) \right|\right)\left(\sum_{\epsilon \in E}\left|\e\left[\log\left(2Y_\eps\right)|\Amostra\right]\right|\right) $$

$\particao$ has VC-dimension $2$, therefore by generalized versions of Glivenko-Cantelli Theorem, the first term is $O(1/\sqrt{n})$ in probability. By Lemma \ref{lemma:thm_2} $\frac{\left(\sum\limits_{\epsilon \in E}\left|\e\left[\log\left(2Y_\eps\right)|\Amostra\right]\right|\right)}{\sqrt{n}} \rightarrow 0$ in probability. It follows that the product converges to $0$ and the theorem is proved. 

\end{proof}

\begin{proof}[Proof of Lemma \ref{lemma:thm_2}]

Consider $L_n = \frac{\log_2 n}{\delta}$ a fixed truncation level of the partition tree for $2 < \delta < \beta$. Also, let $\bar{\eps} = (\eps_1 \dots \eps_{k-1})$. By Proposition \ref{prop:upper-bound-summand} it holds:

$$\left(\sum_{\epsilon \in E}\ex\left[\left|\e\left[\log\left(2Y_\eps\right)|\Amostra\right]\right|\right]\right) \leq \sum_{\epsilon: \len \leq L_n}C + \sum_{\epsilon: \len > L_n} \e\left(\frac{N_{\bar{\eps}}}{\al}\right) = $$
$$2M (2^{L_n}-1) + \sum_{l > L_n} \left(\frac{2n}{\al}\right) \leq  2M(n^{\frac{1}{\delta}}-1) + \sum_{l > L_n} \left(\frac{2n}{\al}\right) \leq 2C(n^{1/\delta}) + 2n^{1-\frac{\beta}{\delta}}$$

therefore $\left(\sum_{\epsilon \in E}\left|\e\left[\log\left(2Y_\eps\right)|\Amostra\right]\right|\right)n^{-1/2} \rightarrow 0$ in $L_1$ and in probability.

\end{proof}

%$$P_X(\BO)\log\left(\frac{n P_X(\BE) + 2\al}{2(n P_X(\BO) + \al - 2)}\right) + P_X(\BU)\log\left(\frac{n P_X(\BE) + 2\al}{2(n P_X(\BU) + \al - 2)}\right) = $$

% $$\e\left[\sum_{\epsilon \neq \epsilon'} \left(\sum_{i=1}^n \sum_{j=1}^n \mathbb{I}_{B_{\epsilon'}} (X_j) \mathbb{I}_{B_\epsilon}(X_i) \right) \log(2 Y_{\eps}) \log(2 Y_{\eps'})\right] + \e\left[\sum_{\epsilon}  \left(\sum_{i=1}^n \sum_{j=1}^n \mathbb{I}_{B_{\epsilon}}(X_j) \mathbb{I}_{B_\epsilon}(X_i)\right) \log(2 Y_{\eps})^2\right] = $$

% $$\left(\sum_{i=1}^n \sum_{j=1}^n \mathbb{I}_{B_{\epsilon'}} (X_j) \mathbb{I}_{B_\epsilon}(X_i) \right) = \begin{cases}
% \nne^2, & \text{if $\epsilon = \epsilon'$}\\
% N_{\min{\eps, \epsl}} + \nne\nnel, & \text{if $\epsilon \neq \epsilon'$}
% \end{cases}$$

% $$\e\left[\sum_{i=1}^n \sum_{j=1}^n \log \theta(X_i) \log \theta(X_i)\right] = \sum_{\epsilon} \nne^2 \e\left[\log(2 Y_\epsilon)^2 | \Amostra\right] + 2\sum_{\epsilon' < \eps}(N_{\min{\eps, \epsl}} + \nne\nnel)\e\left[\log(2 Y_\epsilon) | \Amostra\right] \e\left[\log(2 Y_{\epsl}) | \Amostra\right] + $$

%%%%%%%%%%%%%%%%%%%%%%%%%%%%%%%%%%%%%%%%%%%%%%
%% Funding information, if any,             %%
%% should be provided in the                %%
%% funding section.                         %%
%%%%%%%%%%%%%%%%%%%%%%%%%%%%%%%%%%%%%%%%%%%%%%
\begin{funding}
Rafael Bassi Stern is grateful for the financial support of CNPq (grant 313557/2025-0) and University of São Paulo (PRPI/USP 58/2023), and produced this work as part of the activities of FAPESP, Brazil Research, Innovation and Dissemination Center for Neuromathematics (grant 2013/07699-0).

\end{funding}

%%  MR numbers will be added by VTeX.                      %%
%%                                                         %%
%%  Use \cite{...} to cite references in text.             %%
%%                                                         %%
%%%%%%%%%%%%%%%%%%%%%%%%%%%%%%%%%%%%%%%%%%%%%%%%%%%%%%%%%%%%%

%% if your bibliography is in bibtex format, uncomment commands:
\bibliographystyle{imsart-nameyear} % Style BST file (imsart-number.bst or imsart-nameyear.bst)
\bibliography{bibtex}       % Bibliography file (usually '*.bib')

\end{document}